\documentclass{commat}

\usepackage{graphicx}

\title[On Lie algebras associated with a spray]{On Lie algebras associated with a spray: Case of infinitesimal isometries of a Riemannian manifold}

\author{Manelo Anona and Hasina Ratovoarimanana}

\affiliation{
    \address{Manelo Anona --
Department of Mathematics and Computer Science,
Faculty of Science, University of Antananarivo,
Antananarivo 101, PB: 906, MADAGASCAR
        }
    \email{mfanona@yahoo.fr
    }
    \address{Hasina Ratovoarimanana--Department of Mathematics and Computer Science,
    	Faculty of Science, University of Antananarivo,
    	Antananarivo 101, PB: 906, MADAGASCAR
        }
    \email{hasinaratovo5@gmail.com
    }}

\abstract{
    The Lie algebra of infinitesimal isometries of a Riemannian manifold contains at most two commutative ideals. One coming from the horizontal nullity space of the Nijenhuis tensor of the canonical connection, the other coming from the constant vectors fields  independent of the Riemannian metric.
    }

\keywords{%
Differentiable manifold, Riemannian manifold, Lie algebra, Spray, Nijenhuis tensor, Infinitesimal isometries.
    }

\msc{53XX, 17B66, 53B05, 53C08.
    }

\VOLUME{30}
\NUMBER{2}
\firstpage{13}
\DOI{https://doi.org/10.46298/cm.9007}

\begin{paper}

\section{Introduction}

Let $ M $ be a paracompact differentiable manifold of dimension $ n\geq 2 $ and of class $ \mathcal{C}^\infty $, $ S $ a spray on $ M $. To study the Lie algebra of projectable vector fields which commute with $ S $, which will be denoted by $ \overline{A_S} $ cf.\cite{LOO}, we have associated with $ S $ the  vector $1-$form $ J $ defining the tangent structure on $ M $. The  vector $ 1-$form $ \Gamma = [J, S] $ is considered as a connection in the sense of \cite{GRI}. In \cite{ANO1}, we have shown that the elements of $ \overline{A_S} $, belonging to the horizontal nullity space of the curvature $ R $ of $ \Gamma $, form a commutative ideal of $ \overline{A_S} $. In \cite{ANO2}, we found that some constant elements of $ \overline{A_S} $ can constitute a commutative ideal. In the present study, we show that there are two possible commutative ideals. A Lie algebra of infinitesimal isometries $ \overline{A_g} $ of a Riemannian manifold is semi-simple if and only if the horizontal nullity space of the Nijenhuis tensor of $ \Gamma $ is zero and that the derived ideal  $[\overline{A_g},\overline{A_g}] $ coincides with $ \overline{A_g} $.
To illustrate our results, we give some examples of $ \overline{A_S} $ and $\overline{A_g}$.

\section{Preliminaries}
Let $ K $ and $ L $ be two vector $ 1- $form on a manifold $ M $, $ \chi(M) $ the set of vector fields on $ M $. The bracket $ [K, L] $ cf.\cite{FN1} is written
\begin{eqnarray*}
	[K,L](X,Y)&=&[KX,LY]+[LX,KY]+KL[X,Y]+LK[X,Y]-K[LX,Y]\\&&-L[KX,Y]-K[X,LY]-L[X,KY]
\end{eqnarray*}
for all $X,\ Y\in\chi(M)$.

The bracket $ N_L = \frac{1}{2} [L, L] $ is called the Nijenhuis tensor of $ L $. The Lie derivative of $ X $ applied to $ L $ is written
\begin{equation*}
[X,L]Y=[X,LY]-L[X,Y].
\end{equation*}
The exterior derivation $ d_L $ is defined by $ d_L = [i_L, d] $, where $ i_L $ the inner product with respect to $ L $.

Let $ \Gamma $ be a connection in the sense of \cite{GRI}. The  vector $1-$form $ \Gamma $ is an almost product structure $ (\Gamma^2 = I$, $ I $ being the identity vector $ 1-$form). Noting
\begin{eqnarray*}
	h=\frac{1}{2}(I+\Gamma)\ \text{and} \ v=\frac{1}{2}(I-\Gamma), 
\end{eqnarray*}
The vector $ 1-$form $ h $ is the horizontal projector corresponding to the eigenvalue $ + 1 $, and $ v $ the vertical projector for the eigenvalue $ -1 $. The curvature of $ \Gamma $ is defined by $ R = \frac{1}{2}[h, h] $ which is also equal to $ \frac{1}{8}[\Gamma, \Gamma] $.

The Lie algebra $ A_\Gamma $ is defined by \begin{eqnarray*}
	A_\Gamma=\{X\in\chi(TM)\ \text{such that}\ [X,\Gamma]=0 \}.
\end{eqnarray*}

The nullity space of the curvature $ R $ is:\begin{eqnarray*}
	\mathcal{N}_R=\{X\in \chi(TM)\ \text{ such that}\ R(X,Y)=0,\ \forall\ Y\in \chi(TM)\}.
\end{eqnarray*}
In local natural coordinates on an open set $ U $ of $ M $, $(x^i,y^j)$ are the coordinates on $ TU $, a spray $ S $ is written

\begin{equation*}
S= y^i\frac{\partial}{\partial x^i}-2G^i(x^1,\ldots,x^n,y^1,\ldots,y^n)\frac{\partial}{\partial y^i}.
\end{equation*}
Let $ A_S = \{X \in \chi(TM) \ \text{such that} \ [X, S] = 0 \} $. The equation $ [X, S] = 0 $ implies that the projectable elements of $ A_S $ are of the form
\begin{equation*}
\overline{X}=X^i(x)\frac{\partial}{\partial x^i}+y^j \frac{\partial X^i(x)}{\partial x^j}\frac{\partial}{\partial y^i}.
\end{equation*}
Let $ \overline{\chi (M)} $ denote the complete lift on the tangent bundle $ TM $ of $ \chi(M) $ on $ M $. The projectable elements of $ A_S $ are in $ \overline{\chi(M)} $.
According to Jacobi's identity cf.\cite {FN1}, we can write with the vector $ 1- $form $ J $ defining the tangent structure of $ M $,
\begin{equation*}
[[\overline{X},S],J]+[[S,J],\overline{X}]+[[J,\overline{X}], S]=0.
\end{equation*}
By hypothesis we have $ [\overline{X}, S] = 0 $ and, $ [J, \overline{X}] = 0 $ according to a result of \cite{LEH}, we obtain
\begin{equation*}
[\overline{X},\Gamma]=0,
\end{equation*}
with $\Gamma=[J,S]$. Let $ C $ denote the Liouville field on the tangent bundle $ TM $, the homogeneity of $ \overline{X} $ ($ [C, \overline{X}] = 0 $) leads us to study $ \overline{A_S} $ taking $ S $ such that $ [C, S] = S $, taking into account $ [C, J] = - J $.\\
Thus the connection $ \Gamma $ is linear and without torsion according to \cite{GRI}. For a  connection $ \Gamma = [J, S] $, the coefficients of $ \Gamma $ become $ \Gamma^j_i = \frac{\partial G^j} {\partial y^i} $ and the projectors horizontal and vertical are
\begin{eqnarray*}
	\begin{cases}
		h(\frac{\partial}{\partial x^i})=\frac{\partial}{\partial x^i}-\Gamma^j_i\frac{\partial}{\partial y^j}\\
		h(\frac{\partial}{\partial y^j})=0
	\end{cases}
	\begin{cases}
		v(\frac{\partial}{\partial x^i})=\Gamma^j_i\frac{\partial}{\partial y^j}\\
		v(\frac{\partial}{\partial y^j})= \frac{\partial}{\partial y^j}
	\end{cases}
	i,j\in\{1,\ldots,n\}.\end{eqnarray*}

The curvature $ R = \frac{1}{2}[h, h] $ is then
\[
    R = \frac{1}{2}R^k_{ij}dx^i \wedge dx^j \otimes \frac{\partial}{\partial y^k}
    \quad \textup{ where } \quad
    R^k_{ij} = \frac{\partial \Gamma^k_i}{\partial x^j} - \frac{\partial \Gamma^k_j}{\partial x^i} + \Gamma^l_i \frac{\partial \Gamma^k_j}{\partial y^l} - \Gamma^l_j \frac{\partial \Gamma^k_i}{\partial y^l},
\]
for each $i,j,k,l \in \{1, \dotsc, n \}$
As the functions $ G^k $ are homogeneous of degree $ 2 $, the coefficients $ \Gamma^k_{ij} = \frac{\partial^2 G^k} {\partial y^i \partial y^j} $ do not depend on $ y^i $, $ i \in \{1, \ldots, n \} $. We then have $ R^ k_{ij} = y^l R^k_ {l, ij} (x) $, the $ R^k_{l, ij} (x) $ depend only on the coordinates of the manifold $ M $.

\begin{proposition}[\cite{ANO1}]\label{P.1}
	The Lie algebra $ \overline{A_ {S}} $ coincides with $ \overline{A_{\Gamma}} $.
\end{proposition}
\begin{proof}
	See Proposition 9 of \cite{ANO1}.
\end{proof}

\section{The horizontal elements of $ \overline{A_\Gamma} $}
\begin{proposition}\label{P3.1}
	The elements of $ A_\Gamma $ are projectable vector fields.
\end{proposition}
\begin{proof}
	A vector field $ X \in A_\Gamma $ means $ [X, \Gamma] = 0 $. By definition, the horizontal projector $ h $ is $ h = \frac{I + \Gamma} {2} $. The relation $ [X, \Gamma] = 0 $ is equivalent to $ [X, h] = 0 $. \\
	By expanding $ [X, h] = 0 $, we have for all $ Y \in \chi(TM) $
	\begin{equation*}
	[X,hY]=h[X,Y].
	\end{equation*}
	If $ Y $ is a vertical vector field, we find $ hY = 0 $. the above relation becomes $ h[X, Y] = 0 $ for any vertical vector field, i.e. $ X $ is a projectable vector field.
\end{proof}
We will denote in the following by $ H^\circ $ the set of horizontal and projectable vector fields.
\begin{proposition}[\cite{ANO1}]\label{P3.2}
	Let $ X $ be a projectable vector field. The following two relationships are equivalent:
	\begin{enumerate}
		\item[$(i)$] $[hX,J]=0$;
		\item[$(ii)$] $[JX,h]=0$.
	\end{enumerate}
\end{proposition}
\begin{proposition}\label{P3.3}
	Let $A^h_\Gamma=A_\Gamma\cap H^\circ$ and $\overline{A_\Gamma}^h=A^h_\Gamma\cap\overline{\chi(M)}$. The horizontal vector fields $ \overline{A_\Gamma}^h $ of $ \overline{A_\Gamma} $ form a commutative ideal of $ \overline{A_\Gamma} $. The dimension of $ \overline{A_\Gamma}^h $ corresponds to the dimension of $ A^h_\Gamma $ if the rank of $ A^h_\Gamma $ is constant.
\end{proposition}
\begin{proof}
	From the result of \cite{RRA1}, $ A^h_\Gamma $ is an ideal of $ A_\Gamma $, so $ \overline{A_\Gamma}^h = A^h_\Gamma \cap \overline{\chi(M)} $ is an ideal of $ A_\Gamma \cap \overline{\chi(M)} = \overline{A_\Gamma} $. On the other hand, we have $ v [hZ, hT] = R (Z, T) $, $ \forall \ Z, T \in \chi(TM) $. We then obtain $ v[\overline {X}, \overline{Y}] = 0 $ for $\overline{X},\overline{Y}\in\overline{A_\Gamma}^h$.\\
	According to propositions \ref{P3.2} and 2 of \cite{ANO1}, we have $ J [\overline{X}, \overline{Y}] = 0 $ for $ \overline{X}, \overline{Y} \in \overline{A_ \Gamma}^h $, given $ [J, \Gamma] = 0 $. The horizontal and vertical parts of $ [\overline{X}, \overline{Y}] $ are zero, we find $ [\overline{X}, \overline{Y}] = 0 $.\\
	For the existence of such an element of $ \overline{A_\Gamma}^h $, we must have $ h(\overline{X}) = \overline{X} $, that is to say
	\begin{eqnarray*}
		X^i(x)\frac{\partial}{\partial x^i} -X^i(x)\Gamma_i^j\frac{\partial}{\partial y^j}=X^i(x)\frac{\partial}{\partial x^i}+y^l\frac{\partial X^j(x)}{\partial x^l}\frac{\partial}{\partial y^j}.
	\end{eqnarray*}
	The system of equations to be solved becomes
	\begin{eqnarray}\label{3.3.1}
	\frac{\partial X^j(x)}{\partial x^l}=-X^i(x)\Gamma_{il}^j,\ \text{with}\  i,\ j,\ l\in \{1,\ldots,n\}.
	\end{eqnarray}
	The compatibility condition of such a system of equations according to the Frobenius theorem is
	\begin{eqnarray*}
		X^l(\frac{\partial \Gamma_{li}^k}{\partial x^j}-\frac{\partial \Gamma_{lj}^k}{\partial x^i}+\Gamma_{li}^s\Gamma_{sj}^k-\Gamma_{lj}^s\Gamma_{si}^k)=0, \  i,\ j,\, k,\ l,\ s\in \{1,\ldots,n\}.
	\end{eqnarray*}
	that is, $ X^l R^k_{l, ij} = $ 0. This condition is satisfied if $ X \in H^\circ \cap \mathcal{N}_R = A^h_\Gamma $. The ideal of $ A^h_\Gamma $ is a module over the smooth functions of $ M $ and involutive. On the integral sub-manifold defined by $ A^h_\Gamma $, the system (\ref{3.3.1}) admits solutions in the number of the dimension of the said sub-manifold.
\end{proof}

\section{Lie algebras of infinitesimal isometries}
Let $ E $ be an energy function, a function $\mathcal{T}M=TM-\{0\}$ into $\mathbb{R}^+$, with $E(0)=0$, of class $\mathcal{C}^\infty$ on $\mathcal{T}M$, of class $\mathcal{C}^2$ on the null section, and homogeneous of degree two such that the $ 2-$ form $ \Omega = dd_JE $ being of maximum rank. The canonical spray $ S $ is defined by
\begin{equation*}
i_Sdd_JE=-dE,
\end{equation*}
the derivation $ i_S $ is the inner product with respect to $ S $. The $ 2- $ form $ \Omega $ defines a Riemannian metric $g$ on the vertical bundle:
$$g(JX,JY)=\Omega(JX,Y)$$ for all $X$, $Y\in\chi(TM)$. With a natural local coordinate system on an open set $ U $ of $ M $, $ (x^i, y^j) $ are the system of coordinates on $ TU $, the function $ E $ is written
\begin{equation*}
E=\frac{1}{2}g_{ij}(x^1,\ldots,x^n)y^iy^j,
\end{equation*}
where $ g_{ij} (x^1, \ldots, x^n) $ are positive functions such that the symmetric matrix whose $(i,j)$-th entry is $g_ {ij} (x^1, \ldots, x^n)$ is invertible. The relation $ i_Sdd_JE = -dE $ gives the spray $ S $
\begin{equation*}
S=y^i\frac{\partial}{\partial x^i}-2G^i(x^1,\ldots,x^n,y^1,\ldots,y^n)\frac{\partial}{\partial y^i}
\end{equation*}
with
\begin{equation*}
G_k=\frac{1}{2}y^i y^j\gamma_{ikj}
\ \text{where}\ 
\gamma_{ikj}=\frac{1}{2}(\frac{\partial g_{kj}}{\partial x^i}+\frac{\partial g_{ik}}{\partial x^j}-\frac{\partial g_{ij}}{\partial x^k}),
\end{equation*}
by
\begin{equation*}
\gamma^k_{ij}=g^{kl}\gamma_{ilj},
\end{equation*}
we have 
\begin{equation*}
G^k=\frac{1}{2}y^i y^j\gamma^k_{ij},\ i,j,k\in\{1,\ldots,n\}.
\end{equation*}
\begin{definition}
	A vector field $ X $ on a Riemannian manifold $ (M, E) $ is called an infinitesimal automorphism of the symplectic form $ \Omega $ if $ L_X \Omega = 0 $, where $ L_X $ is the Lie derivative  with respect to $ X $.
\end{definition}
We notice that the canonical spray $ S $ of $ (M, E) $ is an infinitesimal automorphism of the symplectic form $ \Omega $. The set of infinitesimal automorphisms of $ \Omega $ forms a Lie algebra. We denote this Lie algebra by $ A_g $, even formed by the projectable vector fields, it is in general of infinite dimension cf. \cite{ANO1}.
\begin{proposition}[\cite{ANO1}]\label{P4.2} Let $\overline{A_{g}}$ defined as $A_{g}\cap\overline{\chi(M)}$, we have 
	\begin{enumerate}
		\item[a)] $X\in \overline{A_{g}}$ if and only if $ X $ is a projectable vector field such that  $X\in A_{g}$ and $L_{X}E=0$;
		\item[b)] $\overline{A_{g}}\subset \overline{A_\Gamma}$; the horizontal elements of $ \overline{A_\Gamma} $ form a commutative ideal of $ \overline{A_{g}} $.
		\item[c)] The elements of $ \overline{A_{g}} $ are Killing fields of the projectable vectors of the metric $ g $ belonging to $ \overline{A_\Gamma} $. The dimension of $ \overline{A_g} $ is at most equal to $ \frac{n (n + 1)}{2} $. 
	\end{enumerate}
\end{proposition}
\begin{proposition}\label{P4.3}
	On a Riemannian manifold $ (M, E) $, the horizontal nullity space of the curvature $ R $ is generated as a module by the projectable vector fields belonging to this nullity space and, orthogonal to the image space $ Im R $ of the curvature $ R $.
\end{proposition}
\begin{proof}
	If the potential $ R^\circ = i_S R $ is zero, then the curvature $ R $ is zero, in this case the horizontal space $ Im h $ is the horizontal nullity space of the curvature $ R $, isomorphic to $ \chi(U) $, $ U $ being an open set of $ M $.
 
	In the following, we assume that $ R^\circ \neq 0 $. The vertical vector field $ JX $ is orthogonal to the image $ Im R $ of the curvature if and and only if the curvature of the connection $ D $ of Cartan $ \mathfrak{R}(S, X) Y = 0 $ $ \forall Y \in \chi(TM) $ cf.\cite{ANO1}. We obtain $R(X,Y)=R^\circ[JY,X]$ $\forall Y\in \chi(TM)$. As $ R $ is a semi-basic vector $ 2-$form, the above relation is only possible if $ X = S $ or $ X \in hN_R $, then $ X $ is generated as a module by projectable vector fields  in $ hN_R $.
\end{proof}
\begin{theorem}\label{T4.4}
	The Lie algebra $ \overline{A_g} $ is semi-simple if and only if the horizontal nullity space of the Nijenhuis tensor of $ \Gamma $ is zero and, the derived ideal  $[\overline{A_g},\overline{A_g}] $ coincides with $ \overline{A_g} $.
\end{theorem}
\begin{proof}
	If the Lie algebra $ \overline{A_g} $ is semi-simple, any commutative ideal of $ \overline{A_g} $ is zero by definition. According to propositions \ref{P3.3} and \ref{P4.3} the horizontal nullity space of the Nijenhuis tensor of $ \Gamma $ is zero. The derived ideal $[\overline{A_g},\overline{A_g}]$ coincides with $ \overline{A_g} $ according to a classical result.\\
	Conversely, if $ \overline{X} \in \overline{A_\Gamma} $, we have $[\overline{X},h]=0$. According to the identity of Jacobi cf.\cite{FN1} $ [\overline{X}, [h, h]] = 0 $, that is to say, $ [\overline{X}, R ] = 0 $, we then have $ [\overline{X}, R(Y, Z)] = R ([\overline{X}, Y], Z) + R(Y, [\overline {X}, Z]) $, for all $ Y, Z \in \chi(TM) $. If $ \overline {X} $ and $ \overline{Y} $ are elements of a commutative ideal of $ \overline{A_\Gamma} $, we find 
	\begin{equation}\label{4.4.1}
	[\overline{X},R(\overline{Y},Z)]=R(\overline{Y},[\overline{X},Z]),\ \forall Z\in\chi(TM).
	\end{equation}
	If the horizontal nullity space of the curvature $ R $ is zero, the  semi-basic vector $ 2-$form $ R $ is non-degenerate. The only possible case for equation (\ref{4.4.1}) is that the commutative ideal of $ \overline{A_g} $ is at most formed by constant vector fields $ \frac{\partial} {\partial x^i} $, $ i \in \{1, \dots, n \} $ such that $ \frac{\partial G^k} {\partial x^i} = 0 $ for all $ k \in \{1, \dots, n \} $ from \cite{ANO2}. These constant vector fields can only form an ideal of affine vector fields independent of the other elements of $ \overline{A_g} $. The derived ideal  $ [\overline{A_g},\overline{A_g}]$ never coincides with $ \overline{A_g} $.
	Hence the result.
\end{proof}
\begin{example}
	We take $ M = \mathbb{R}^3 $ and the energy function is written:
	$$E=\frac{1}{2}(e^{x^3}(y^1)^2+e^{x^3}(y^2)^2+(y^3)^2).$$
	The canonical spray of $ E $ is written:
	\begin{equation*}
	S= y^1\frac{\partial}{\partial x^1}+y^2\frac{\partial}{\partial x^2}+y^3\frac{\partial}{\partial x^3} - y^1y^3\frac{\partial}{\partial y^1}- y^2y^3\frac{\partial}{\partial y^2}+\frac{e^{x^3}}{2}((y^1)^2+(y^2)^2)\frac{\partial}{\partial y^3}.
	\end{equation*}
	The non-zero coefficients of $ \Gamma $ are
	\begin{eqnarray*}
		&&\Gamma^1_1=\frac{y^3}{2},\ \Gamma^1_3= \frac{y^1}{2},\ \Gamma^2_2=\frac{y^3}{2},\ \Gamma^2_3=\frac{y^2}{2},\\  &&\Gamma^3_1=-\frac{e^{x^3}y^1}{2},\ \Gamma^3_2=- \frac{e^{x^3}y^2}{2}.
	\end{eqnarray*}
	The horizontal fields are generated as a module by
	\begin{eqnarray*}
		&&\frac{\partial}{\partial x^1}-\frac{y^3}{2}\frac{\partial}{\partial y^1}+\frac{e^{x^3}y^1}{2}\frac{\partial}{\partial y^3},\\
		&&\frac{\partial}{\partial x^2}-\frac{y^3}{2}\frac{\partial}{\partial y^2}-\frac{e^{x^3}y^2}{2}\frac{\partial}{\partial y^3},\\
		&&\frac{\partial}{\partial x^3}-\frac{y^1}{2}\frac{\partial}{\partial y^1}-\frac{y^2}{2}\frac{\partial}{\partial y^3}.
	\end{eqnarray*}
	The horizontal nullity space of the curvature is zero.
	The Lie algebra $ \overline {A_\Gamma} $ is generated as Lie algebra by:
	\begin{eqnarray*}
		e_1&=& -\frac{x^2x^1}{2}\frac{\partial}{\partial x^1}+(e^{-x^3}+\frac{(x^1)^2}{4}-\frac{(x^2)^2}{4})\frac{\partial}{\partial x^2}+x^2\frac{\partial}{\partial x^3}-\frac{y^2 x^1+y^1 x^2}{2}\frac{\partial}{\partial y^1}\\&&-\frac{(-x^1y^1+x^2y^2+2y^3 e^{-x^3})}{2}\frac{\partial}{\partial y^2}+y^2\frac{\partial}{\partial y^3},\\
		e_2&=& (-2e^{-x^3}+\frac{(x^1)^2}{2}-\frac{(x^2)^2}{2})\frac{\partial}{\partial x^1}+x^1x^2\frac{\partial}{\partial x^2}-2x^1\frac{\partial}{\partial x^3}\\&&+(-2y^3e^{-x^3}+x^1y^1-x^2y^2)\frac{\partial}{\partial y^1}+(x^1y^2+x^2y^1)\frac{\partial}{\partial y^2}-2y^1\frac{\partial}{\partial y^3},\\
		e_3&=&-x^2\frac{\partial}{\partial x^1}+x^1\frac{\partial}{\partial x^2}-y^2\frac{\partial}{\partial y^1}+y^1\frac{\partial}{\partial y^2},\\
		e_4&=&x^1\frac{\partial}{\partial x^1}+x^2\frac{\partial}{\partial x^2}-2\frac{\partial}{\partial x^3}+y^1\frac{\partial}{\partial y^1}+y^2\frac{\partial}{\partial y^2},\\
		e_5&=&\frac{\partial}{\partial x^2},\\
		e_6&=&\frac{\partial}{\partial x^1}.
	\end{eqnarray*}

        \begin{table}\centering
	\begin{tabular}{|c|c|c|c|c|c|c|}
		\hline
		$[\cdot , \cdot]$& $e_1$& $e_2$& $e_3$& $e_4$& $e_5$& $e_6$\\
		\hline
		$e_1$& $0$& $0$&$\frac{e_2}{2}$&$-e_1$&$\frac{e_4}{2}$&$-\frac{e_3}{2}$\\
		\hline
		$e_2$& $0$& $0$&$-2e_1$&$-e_2$&$-e_3$&$-e_4$\\
		\hline
		$e_3$& $-\frac{e_2}{2}$& $2e_1$&$0$&$0$&$e_6$&$-e_5$\\
		\hline
		$e_4$& $e_1$& $e_2$&$0$&$0$&$-e_5$&$-e_6$\\
		\hline
		$e_5$& $-\frac{e_4}{2}$& $e_3$&$-e_6$&$e_5$&$0$&$0$\\
		\hline
		$e_6$& $\frac{e_3}{2}$& $e_4$&$e_5$&$e_6$&$0$&$0$\\
		\hline
	\end{tabular}
        \caption{Multiplication table of $\overline{A_\Gamma}$}
        \end{table}
	
	The Lie algebra $\overline{A_\Gamma}=\overline{A_g}$ is simple.
\end{example}
\begin{example} \label{eg:R4}
	We take $ M = \mathbb{R}^4 $ and the energy function is written:
	$$E=\frac{1}{2}(e^{x^2}(y^1)^2+(y^2)^2+e^{x^4}(y^3)^2+(y^4)^2).$$
	The canonical spray of $ E $ is written:
	\begin{equation*}
	S= y^1\frac{\partial}{\partial x^1}+y^2\frac{\partial}{\partial x^2}+y^3\frac{\partial}{\partial x^3}+y^4\frac{\partial}{\partial x^4} - y^1y^2\frac{\partial}{\partial y^1}+\frac{e^{x^2}}{2} (y^1)^2\frac{\partial}{\partial y^2}- y^3y^4\frac{\partial}{\partial y^3}+\frac{e^{x^4}}{2} (y^3)^2\frac{\partial}{\partial y^4}.
	\end{equation*}
	The non-zero coefficients of $ \Gamma $ are
	\begin{eqnarray*}
		\Gamma^1_1=\frac{y^2}{2},\ \Gamma^1_2= \frac{y^1}{2},\ \Gamma^2_1=-\frac{e^{x^2}y^1}{2},\  \Gamma^3_2=\frac{y^4}{2},\ \Gamma^3_3=\frac{y^3}{2},\ \Gamma^4_3=- \frac{e^{x^4}y^3}{2}.
	\end{eqnarray*}
	The horizontal fields are generated as a module by
	\begin{eqnarray*}
		&&\frac{\partial}{\partial x^1}-\frac{y^2}{2}\frac{\partial}{\partial y^1}+\frac{e^{x^2}y^1}{2}\frac{\partial}{\partial y^2},\\
		&&\frac{\partial}{\partial x^2}-\frac{y^1}{2}\frac{\partial}{\partial y^1},\\
		&&\frac{\partial}{\partial x^3}-\frac{y^4}{2}\frac{\partial}{\partial y^3}+\frac{e^{x^4}y^3}{2}\frac{\partial}{\partial y^4}\\
		&&\frac{\partial}{\partial x^4}-\frac{y^3}{2}\frac{\partial}{\partial y^3},\\.
	\end{eqnarray*}
	The horizontal nullity space of the curvature is zero.
	The Lie algebra $ \overline {A_\Gamma} $ is generated as Lie algebra by:
	\begin{eqnarray*}
		e_1&=& -(-e^{-x^2}+\frac{(x^1)^2}{4})\frac{\partial}{\partial x^1}+x^1\frac{\partial}{\partial x^2}-(\frac{-x^1y^1}{2}+y^2 e^{-x^2})\frac{\partial}{\partial y^1}+y^1\frac{\partial}{\partial y^2},\\
		e_2&=& -\frac{x^1}{2}\frac{\partial}{\partial x^1}+\frac{\partial}{\partial x^2}-\frac{y^2}{2}\frac{\partial}{\partial y^1},\
		e_3=\frac{\partial}{\partial x^1},\\
		e_4&=& -(-e^{-x^4}+\frac{(x^3)^2}{4})\frac{\partial}{\partial x^3}+x^3\frac{\partial}{\partial x^4}-(\frac{-x^3y^3}{2}+y^4 e^{-x^4})\frac{\partial}{\partial y^3}+y^3\frac{\partial}{\partial y^4},\\
		e_5&=& -\frac{x^3}{2}\frac{\partial}{\partial x^3}+\frac{\partial}{\partial x^4}-\frac{y^3}{2}\frac{\partial}{\partial y^3},\
		e_6=\frac{\partial}{\partial x^3},\\
	\end{eqnarray*}
	
	\begin{table}\centering
	\begin{tabular}{|c|c|c|c|c|c|c|}
		\hline
		$[\cdot,\cdot]$& $e_1$& $e_2$& $e_3$& $e_4$& $e_5$& $e_6$\\
		\hline
		$e_1$& $0$& $-e_1$&$\frac{e_2}{2}$&$0$&$0$&$0$\\
		\hline
		$e_2$& $e_1$& $0$&$-e_3$&$0$&$0$&$0$\\
		\hline
		$e_3$& $-\frac{e_2}{2}$& $e_3$&$0$&$0$&$0$&$0$\\
		\hline
		$e_4$& $0$& $0$&$0$&$0$&$-e_4$&$\frac{e_5}{2}$\\
		\hline
		$e_5$& $0$& $0$&$0$&$e_4$&$0$&$-e_6$\\
		\hline
		$e_6$& $0$& $0$&$0$&$-\frac{e_5}{2}$&$e_6$&$0$\\
		\hline
	\end{tabular}
        \caption{Multiplication table of $\overline{A_\Gamma}$ in Example~\ref{eg:R4}}
        \end{table}
	$$\overline{A_S}=\overline{A_S}_s^1\oplus\overline{A_S}_s^2,$$
	such that  $\overline{A_S}_s^1=\{e_{1},e_{2},e_{3}\}\cong sl(2)$ and $\overline{A_S}_s^2=\{e_{4},e_{5},e_{6}\}\cong sl(2)$ are simple.
	The Lie algebra $\overline{A_\Gamma}=\overline{A_g}$ is semi-simple.
\end{example}

\section{An example of a Lie algebra $ \overline{A_S} $ of maximal rank $ n^2 + n $ with a Lie algebra $\overline{A_g}$ of maximal rank $\frac{n^2 + n}{2}$ of different natures.}
We take $ M = \mathbb {R} ^ 3 $ and the energy function is 
\begin{equation*}
E=\frac{1}{2}(e^{x^1}(y^1)^2+e^{x^2}(y^2)^2+e^{x^3}(y^3)^2)
\end{equation*}
The canonical spray of $ E $ is written:
\begin{equation*}
S= y^1\frac{\partial}{\partial x^1}+y^2\frac{\partial}{\partial x^2}+y^3\frac{\partial}{\partial x^3} - \frac{(y^1)^2}{2}\frac{\partial}{\partial y^1}--\frac{(y^2)^2}{2}\frac{\partial}{\partial y^2}-\frac{(y^3)^2}{2}\frac{\partial}{\partial y^3}.
\end{equation*}
The non-zero coefficients  $ \Gamma^j_i $ of $ \Gamma $ are
\begin{equation*}
\Gamma^1_1= \frac{y^1}{2},\ \Gamma^2_2 = \frac{y^2}{2},\ \Gamma^3_3 = \frac{y^3}{2}.
\end{equation*}
The basis of the horizontal space of $ \Gamma $ is written
\begin{eqnarray*}
	\frac{\partial}{\partial x^1}-\frac{y^1}{2}\frac{\partial}{\partial y^1},\ \frac{\partial}{\partial x^2}-\frac{y^2}{2}\frac{\partial}{\partial y^2},\ \frac{\partial}{\partial x^3}-\frac{y^3}{2}\frac{\partial}{\partial y^3}.
\end{eqnarray*}
The curvature $ R $ is zero.
The Lie algebra $ \overline{A_S} $ is generated as Lie algebra by:
\begin{eqnarray*}
	&&e_1=\frac{\partial}{\partial x^1},\
	e_2= e^{\frac{x^2-x^1}{2}}(\frac{\partial}{\partial x^1}-\frac{(y^1-y^2)}{2}\frac{\partial}{\partial y^1}),\
	e_3=e^{\frac{-x^1}{2}}(\frac{\partial}{\partial x^1}-\frac{y^1}{2}\frac{\partial}{\partial y^1}),\\
	&&e_4=e^{\frac{x^3-x^1}{2}}(\frac{\partial}{\partial x^1}-\frac{(y^1-y^3)}{2}\frac{\partial}{\partial y^1}),\
	e_5=e^{\frac{x^1-x^2}{2}}(\frac{\partial}{\partial x^2}-\frac{(y^2-y^1)}{2}\frac{\partial}{\partial y^2}),\\
	&&e_6=\frac{\partial}{\partial x^2},\
	e_7=e^{\frac{-x^2}{2}}(\frac{\partial}{\partial x^2}-\frac{y^2}{2}\frac{\partial}{\partial y^2}),\ e_8=e^{\frac{x^3-x^2}{2}}(\frac{\partial}{\partial x^2}-\frac{(y^2-y^3)}{2}\frac{\partial}{\partial y^2}),\\ &&e_9=e^{\frac{x^1-x^3}{2}}(\frac{\partial}{\partial x^3}-\frac{(y^3-y^1)}{2}\frac{\partial}{\partial y^3}),\ e_{10}=e^{\frac{x^2-x^3}{2}}(\frac{\partial}{\partial x^3}-\frac{(y^3-y^2)}{2}\frac{\partial}{\partial y^3}),\\
	&&e_{11}=\frac{\partial}{\partial x^3},\
	e_{12}=e^{\frac{-x^3}{2}}(\frac{\partial}{\partial x^3}-\frac{y^3}{2}\frac{\partial}{\partial y^3}),\\
\end{eqnarray*}
        \begin{table}\centering
        \resizebox{\textwidth}{!}{%
	\begin{tabular}{|c|c|c|c|c|c|c|c|c|c|c|c|c|}
		\hline
		$[\cdot,\cdot]$& $e_1$& $e_2$& $e_3$& $e_4$& $e_5$& $e_6$& $e_7$& $e_8$& $e_9$& $e_{10}$& $e_{11}$& $e_{12}$\\
		\hline
		$e_1$& $0$& $-\frac{e_2}{2}$&$-\frac{e_3}{2}$&$-\frac{e_4}{2}$&$\frac{e_5}{2}$&$0$& $0$& $0$&$\frac{e_9}{2}$&$0$&$0$&$0$\\
		\hline
		$e_2$& $\frac{e_2}{2}$& $0$&$0$&$0$&$-\frac{e_1}{2}+\frac{e_6}{2}$&$-\frac{e_2}{2}$& $-\frac{e_3}{2}$& $-\frac{e_2}{2}$&$\frac{e_{10}}{2}$&$0$&$0$&$0$\\
		\hline
		$e_3$& $\frac{e_3}{2}$& $0$&$0$&$0$&$\frac{e_7}{2}$&$0$& $0$& $0$&$\frac{e_{12}}{2}$&$0$&$0$&$0$\\
		\hline
		$e_4$& $\frac{e_4}{2}$& $0$&$0$&$0$&$\frac{e_8}{2}$&$0$& $0$& $0$&$-\frac{e_1}{2}+\frac{e_{11}}{2}$&$-\frac{e_2}{2}$&$-\frac{e_4}{2}$&$-\frac{e_3}{2}$\\
		\hline
		$e_5$& $-\frac{e_5}{2}$& $\frac{e_1}{2}-\frac{e_6}{2}$&$-\frac{e_7}{2}$&$-\frac{e_8}{2}$&$0$&$\frac{e_5}{2}$& $0$& $0$&$0$&$\frac{e_9}{2}$&$0$&$0$\\
		\hline
		$e_6$& $0$& $\frac{e_2}{2}$&$0$&$0$&$-\frac{e_5}{2}$&$0$& $-\frac{e_7}{2}$& $-\frac{e_8}{2}$&$0$&$\frac{e_{10}}{2}$&$0$&$0$\\
		\hline
		$e_7$& $0$& $\frac{e_3}{2}$&$0$&$0$&$0$&$\frac{e_7}{2}$& $0$& $0$&$0$&$\frac{e_{12}}{2}$&$0$&$0$\\
		\hline
		$e_8$& $0$& $\frac{e_2}{2}$&$0$&$0$&$0$&$\frac{e_8}{2}$& $0$& $0$&$-\frac{e_5}{2}$&$-\frac{e_6}{2}+\frac{e_{11}}{2}$&$-\frac{e_8}{2}$&$-\frac{e_7}{2}$\\
		\hline
		$e_9$& $-\frac{e_9}{2}$& $-\frac{e_{10}}{2}$&$-\frac{e_{12}}{2}$&$\frac{e_1}{2}-\frac{e_{11}}{2}$&$0$&$0$& $0$& $\frac{e_{5}}{2}$&$0$&$0$&$\frac{e_9}{2}$&$0$\\
		\hline
		$e_{10}$& $0$& $0$&$0$&$\frac{e_{2}}{2}$&$-\frac{e_{9}}{2}$&$-\frac{e_{10}}{2}$& $\frac{e_{12}}{2}$& $\frac{e_{6}}{2}-\frac{e_{11}}{2}$&$0$&$0$&$\frac{e_{10}}{2}$&$0$\\
		\hline
		$e_{11}$& $0$& $0$&$0$&$\frac{e_{4}}{2}$&$0$&$0$& $0$& $\frac{e_{8}}{2}$&$-\frac{e_{9}}{2}$&$-\frac{e_{10}}{2}$&$0$&$-\frac{e_{12}}{2}$\\
		\hline
		$e_{12}$& $0$& $0$&$0$&$\frac{e_{3}}{2}$&$0$&$0$& $0$& $\frac{e_{7}}{2}$&$0$&$0$&$\frac{e_{12}}{2}$&$0$\\
		\hline
	\end{tabular}%
        }
        \caption{Multiplication table of $\overline{A_S}$}
        \end{table}

We then have 
$$ [\overline{A_S}, \overline{A_S}] \neq \overline{A_S}. $$
All derivations of $ \overline{A_S} $ are inner and, $\{e_3, e_7, e_{12}\}$ forms the commutative ideal of $\overline{A_S}$ and include in the horizontal nullity space of the curvature.

The Levi decomposition of $\overline{A_S}$ is $$\overline{A_S}=\overline{A_S}_s+\overline{A_S}_r,$$
such that  $\overline{A_S}_s=\{e_{1}-e_{11},e_{2},e_{4},e_{5},e_{6}-e_{11},e_{8}, e_{9},e_{10}\}$ is simple and $\overline{A_S}_r=\{e_1+e_6+e_{11},e_{3},e_{7},e_{12}\}$ is solvable.
Hence, the Lie algebra $ \overline{A_S} $ is not semi-simple.

The Lie algebra $ \overline{A_g} $ is generated as Lie algebra by
\begin{eqnarray*}
	g_1&=&e^{\frac{x^2-x^1}{2}}(\frac{\partial}{\partial x^1}-\frac{(y^1-y^2)}{2}\frac{\partial}{\partial y^1})-e^{\frac{x^1-x^2}{2}}(\frac{\partial}{\partial x^2}-\frac{(y^2-y^1)}{2}\frac{\partial}{\partial y^2}),\\
	g_2&=&e^{\frac{-x^1}{2}}(\frac{\partial}{\partial x^1}-\frac{y^1}{2}\frac{\partial}{\partial y^1}),\\
	g_3&=&e^{\frac{x^3-x^1}{2}}(\frac{\partial}{\partial x^1}-\frac{(y^1-y^3)}{2}\frac{\partial}{\partial y^1})-e^{\frac{x^1-x^3}{2}}(\frac{\partial}{\partial x^3}-\frac{(y^3-y^1)}{2}\frac{\partial}{\partial y^3}),\\
	g_4&=&e^{\frac{-x^2}{2}}(\frac{\partial}{\partial x^2}-\frac{y^2}{2}\frac{\partial}{\partial y^2})\\
	g_5&=&e^{\frac{x^3-x^2}{2}}(\frac{\partial}{\partial x^2}-\frac{(y^2-y^3)}{2}\frac{\partial}{\partial y^2})-e^{\frac{x^2-x^3}{2}}(\frac{\partial}{\partial x^3}-\frac{(y^3-y^2)}{2}\frac{\partial}{\partial y^3}),\\
	g_6&=&e^{\frac{-x^3}{2}}(\frac{\partial}{\partial x^3}-\frac{y^3}{2}\frac{\partial}{\partial y^3})\\
\end{eqnarray*}

\begin{table}\centering
\begin{tabular}{|c|c|c|c|c|c|c|}
	\hline
	$[\cdot,\cdot]$& $g_1$& $g_2$& $g_3$& $g_4$& $g_5$& $g_6$\\
	\hline
	$g_1$& $0$& $\frac{g_4}{2}$&$\frac{g_5}{2}$& $-\frac{g_2}{2}$& $-\frac{g_3}{2}$& $0$\\
	\hline
	$g_2$& $-\frac{g_4}{2}$& $0$&$-\frac{g_6}{2}$& $0$& $0$& $0$\\
	\hline
	$g_3$& $-\frac{g_5}{2}$& $\frac{g_6}{2}$&$0$& $0$& $\frac{g_1}{2}$& $-\frac{g_2}{2}$\\
	\hline
	$g_4$& $\frac{g_2}{2}$& $0$&$0$& $0$& $-\frac{g_6}{2}$& $0$\\
	\hline
	$g_5$& $\frac{g_3}{2}$& $0$&$-\frac{g_1}{2}$& $\frac{g_6}{2}$& $0$& $-\frac{g_4}{2}$\\
	\hline
	$g_6$& $0$& $0$&$\frac{g_2}{2}$& $0$& $\frac{g_4}{2}$& $0$\\
	\hline
\end{tabular}
\caption{Multiplication table of $\overline{A_g}$}
\end{table}

We have
$$[\overline{A_g},\overline{A_g} ] = \overline{A_g}.$$
The derivations of $ \overline{A_g} $ are inner except for the derivation $
\begin{cases}
D(g_2)=g_2\\
D(g_4)=g_4\\
D(g_6)=g_6
\end{cases}
$ which is outer and, $ \{g_2, \ g_4 \, \ g_6 \} $ forms the commutative ideal of $\overline{A_g}$ and include in the horizontal nullity space of the curvature.
$$[\overline{A_g},\overline{A_g} ] = \overline{A_g}.$$
The Levi decomposition of $\overline{A_g}$ is $$\overline{A_g}=\overline{A_g}_s+\overline{A_g}_i,$$
such that  $\overline{A_g}_s=\{g_{1}, g_{3},g_{5}\}\cong so(3)$ is simple  , $\overline{A_g}_i=\{g_{2},g_{4},g_{6}\}$ is solvable.
Hence, the Lie algebra $ \overline{A_g} $ is not semi-simple.


\EditInfo{January 27, 2022}{June 26, 2022}{Friedrich Wagemann}

\end{paper}